\documentclass[10pt, reqno, a4paper, amsthm]{amsart}

\usepackage{url}
\usepackage{amssymb}
\usepackage{amsmath}
\usepackage{amsfonts} 
\usepackage{mathrsfs}
\usepackage{enumerate} 
\usepackage{verbatim}

\def\ThCon#1{\mathop{
\unitlength=1mm
\begin{picture}(6,4)
\put(0,0){$\leftrightarrow$} \put(1.1,1.2){$*$}\put(3.5,-1){{\tiny
$#1$}}
\end{picture}}}

\def\ssr#1{\mathop{
\unitlength=1mm
\begin{picture}(5,4)
\put(0,0){$\rightarrow$} 
\put(3,-0.6){{\tiny $ _#1$}}
\end{picture}}}

\def\RhCon#1{\mathop{
\unitlength=1mm
\begin{picture}(7,4)
\put(0,0){$\longrightarrow$} \put(2,1.2){$*$}\put(5.1,-1){{\tiny
$#1$}}
\end{picture}}}

\def\RhConAS#1{\mathop{
\unitlength=1mm
\begin{picture}(7,4)
\put(0,0){$\longrightarrow$} \put(2,1.6){{\tiny
$+$}}\put(5.1,-1){{\tiny $ _#1$}}
\end{picture}}}

\newcommand{\gh}{\mathcal{H}}
\newcommand{\gr}{\mathcal{R}}
\newcommand{\gl}{\mathcal{L}}

\newcommand{\gd}{\mathcal{D}}
\newcommand{\gj}{\mathcal{J}}

 \theoremstyle{theorem}
 
 \newtheorem{thm}{Theorem}
 
 \newtheorem{cor}{Corollary}
 \newtheorem{prop}{Proposition}
 \newtheorem{lem}{Lemma}

 \theoremstyle{definition}

\newcommand{\lb}{\langle}
\newcommand{\rb}{\rangle}

\usepackage[active]{srcltx}

\usepackage{verbatim}
\usepackage{enumerate}

\begin{document}

\title[Finite Complete Rewriting Systems]{Finite Complete Rewriting Systems \\ for Regular Semigroups}

\maketitle

\begin{center}

    R. GRAY\footnote{
	  Part of this work was done while this author held an EPSRC Postdoctoral Fellowship at the University of St Andrews.}  

    \medskip

    Centro de \'{A}lgebra, da Universidade de Lisboa, \\ Av. Prof. Gama 		Pinto, 2,  1649-003 Lisboa,  \ Portugal.

    \medskip
    
    \texttt{rdgray@fc.ul.pt}

    \bigskip
    
    \bigskip

    A. MALHEIRO\footnote{
  	This work was developed within the projects POCTI-ISFL-1-143 and PTDC/MAT/69514/2006 of CAUL, financed by FCT and FEDER.}

    \medskip

    Centro de \'{A}lgebra, da Universidade de Lisboa, \\ Av. Prof. Gama 		Pinto, 2,  1649-003 Lisboa,  \ Portugal, and \\ 
    Departamento de Matem\'{a}tica, \ Faculdade de Ci\^{e}ncias e Tecnologia, \\ Universidade Nova de Lisboa, \ 2829-516 Caparica, \ Portugal.

    \medskip

    \texttt{malheiro@cii.fc.ul.pt} \\
\end{center}

\begin{abstract}
It is proved that, given a (von Neumann) regular semigroup with finitely many left and right ideals, if every maximal subgroup is presentable by a finite complete rewriting system, then so is the semigroup. To achieve this, the following two results are proved: the property of being defined by a finite complete rewriting system is preserved when taking an ideal extension by a semigroup defined by a finite complete rewriting system; a completely $0$-simple semigroup with finitely many left and right ideals admits a presentation by a finite complete rewriting system provided all of its maximal subgroups do. 

\

\end{abstract}

\section{Introduction}

This short note is designed to answer a problem posed in \cite[Remark and Open Problem 4.5]{Rus99}. In that paper, the author shows that a regular semigroup with finitely many left and right ideals is finitely presented if and only if all its maximal subgroups are finitely presented. Recall that a semigroup $S$ is said to be (von Neumann) regular if for all $x \in S$ there exists $y \in S$ such that $xyx=x$. Regular semigroups are an important class, which includes many natural examples (e.g. the full transformation semigroup and the full matrix semigroup over a field) and as such have received serious attention in the literature; see \cite{Nambooripad1979} for example.  
Roughly speaking, a regular semigroup is a semigroup with ``lots of'' idempotents (formally, every Green's $\gr$-class and every $\gl$-class of $S$ contains an idempotent). Each of these idempotents is the identity of a maximal subgroup of the semigroup. Consequently it is often the case that the behavior of such semigroups is closely linked to the behavior of their maximal subgroups. The above-mentioned result of Ruskuc is an example of this phenomenon, and there are numerous other analogous results to this where finite presentability is replaced by other standard finiteness properties, including: residual finiteness, local finiteness, periodicity and the property of having soluble word problem; see \cite{Gol75,Rus99}. 

In \cite[Remark and Open Problem 4.5]{Rus99} it is asked whether the same can be established for various finiteness properties relating to homology and rewriting systems, such as being presentable by a finite complete rewriting system, having finite derivation type ($\mathrm{FDT}$), finite cohomological dimension, or satisfying the homological finiteness property $\mathrm{FP_n}$ or $\mathrm{FP_{\infty}}$. In \cite{GM09} we positively answered this question for the property $\mathrm{FDT}$. The questions regarding finite cohomological dimension, $\mathrm{FP_n}$ and $\mathrm{FP_{\infty}}$ are discussed in \cite{GrayPride2010} where is it observed that the corresponding result does not hold for any of these properties. 

Here we consider the property of being presentable by a finite complete rewriting system. Recall that a complete rewriting system for a semigroup is a presentation of a particular kind (both noetherian and confluent) which in particular gives a solution to the word problem for the semigroup (see Section~\ref{sec_prelims} for more details). Therefore it is clearly of considerable interest to develop an understanding of which semigroups are definable by finite complete rewriting systems. Here we prove the following:
\begin{thm}\label{smalltobig}
Let $S$ be a regular semigroup with finitely many left and right ideals. If every maximal subgroup of $S$ is defined by a finite complete rewriting system then $S$ is also defined by a finite complete rewriting system.
\end{thm}     
One source of natural examples of regular semigroups satisfying the hypotheses of Theorem~\ref{smalltobig} is given by, so-called, free regular idempotent generated semigroups of finite regular biordered sets; see \cite{Nambooripad1979} and more recently \cite{Brittenham2009,GrayRuskuc2009}.

The converse of Theorem~\ref{smalltobig} remains open, and may be hard in light of the analogous open problem for groups: is the property of being defined by a finite complete rewriting system inherited by subgroups of finite index? See \cite{PW00b}. Related to this, it is still unknown whether the property of being defined by a finite complete rewriting system is preserved when passing to a subsemigroup with finite complement; see~\cite{Wan98b}.    

Interestingly, both for the main result of \cite{GM09}, and for Theorem \ref{smalltobig} above, regularity is a necessary assumption, i.e. the main result of \cite{Rus00} (which generalises the result mentioned above about finite presentability to non-regular semigroups, with maximal subgroup replaced by the more general concept of Sch\"{u}tzenberger group) does not hold either for $\mathrm{FDT}$, or for the property of being presented by a finite complete rewriting system. A counterexample is given in \cite{GMP10}. 

Our approach to the proof of Theorem~\ref{smalltobig} requires us to introduce a little more theory.

Given a semigroup $S$ and an ideal $I$ of $S$ we can define a
semigroup $S/I$, called the {\it Rees quotient} of $S$ by $I$,
where the elements of $S/I$ are the equivalence classes of the
congruence $\rho_I$ defined on $S$ as follows: for $s,t\in S$, we define $s\rho_I
t$ if either $s=t$ or both $s$ and $t$ are in $I$. 
For more background on basic notions of semigroup theory mentioned here the reader is referred to \cite{How95}.   

Let $S$, $T$ and $U$ be semigroups. The semigroup $S$ is said to
be an {\it ideal extension of $T$ by $U$} if $T$ is isomorphic to
an ideal $T'$ of $S$ and the Rees quotient $S/{T'}$ is isomorphic
to $U$. It is known that an ideal extension of a finitely
presented semigroup by another finitely presented semigroup is
finitely presented; see \cite[Proposition~4.4]{Rus99}. For complete rewriting systems we shall prove the following analogous result:

\begin{thm}\label{theoremIdealExt}
Let $S$ be an ideal extension of a semigroup $T$ by a semigroup $U$. If $T$ and $U$ are both defined by  finite complete rewriting systems then
 $S$ is also defined by a finite complete rewriting system.
\end{thm}

Another result is needed to prove Theorem~\ref{smalltobig}. Recall that a semigroup is said to be \emph{completely ($0$-)simple} if it is ($0$-)simple and has ($0$-)minimal left and right ideals. We shall prove the following:

\begin{thm}\label{comp0simpExt}
Let $S$ be a completely $0$-simple semigroup with finitely many left and right ideals. If every maximal subgroup of $S$ is defined by a finite complete rewriting system  then so is $S$.
\end{thm}

We note that all the non-zero maximal subgroups of a completely $0$-simple semigroup are isomorphic to each other. Also, it is an easy consequence of the Rees-Suschkevitz Theorem (see Section~\ref{sec_comp0simple}) that every finitely generated completely $0$-simple semigroup necessarily has finitely many left and right ideals.  

Using Theorem~\ref{comp0simpExt} we recover the following result originally proved in \cite[Theorem 1.3]{Mal08}.
\begin{cor}\label{compsimp}
Let $S$ be a completely simple semigroup with finitely many left and right ideals, and let $G$ be a maximal subgroup of $S$. If $G$ is defined be a finite complete rewriting system then $S$ is defined by a finite complete rewriting system.
\end{cor}
\begin{proof}
By Theorem \ref{comp0simpExt} the completely $0$-simple semigroup $S^0 = S \cup \{ 0 \}$ given by adjoining a zero element to $S$ is defined by a finite complete rewriting system. Note that $S$ is a subsemigroup of $S^0$ whose complement $S^0\setminus S=\{0\}$ is an ideal.  
It then follows from \cite[Theorem~C]{PW00} that $S$ is defined by a finite complete rewriting system. 
\end{proof}

We now show how Theorem~\ref{smalltobig} may be deduced from Theorems \ref{theoremIdealExt} and \ref{comp0simpExt}.

\begin{proof}[Proof of Theorem~\ref{smalltobig}]
From the assumption that $S$ has finitely many left and right ideals it follows (see for example \cite[Chapter~6]{CP67}) that in $S$ we have $\gj = \gd$ (see \cite[Chapter~2]{How95} for the definition of Green's relations) and hence every principal factor of $S$ is completely $0$-simple.
The proof of Theorem~\ref{comp0simpExt} now goes by induction on the number of $(\gj = \gd)$-classes of $S$.
When $S$ has just one $\gj$-class it follows, since $S$ has finitely many $\gr$- and $\gl$-classes, that $S$ is
 isomorphic to a completely simple semigroup with finitely many left and right ideals which has a maximal subgroup defined by a finite complete rewriting system by assumption.
It follows by Corollary~\ref{compsimp} that the completely simple semigroup $S$ is defined by a finite complete rewriting system.

Now suppose that $S$ has at least two $\gj$-classes. Let $J$ be a maximal $\gj$-class (in the natural ordering of $\gj$-classes $J_x \leq_{\gj} J_y \Leftrightarrow S^1 x S^1 \subseteq S^1 y S^1$). Then $T = S \setminus J$ is an ideal of $S$ where $T$ is regular and, since $T$ is a union of $(\gj = \gd)$-classes of $S$ each of which is regular and contains only finitely many $\gr$- and $\gl$-classes, it follows that
$T$ has strictly fewer $(\gj = \gd)$-classes than $S$. So $T$ is a regular semigroup with finitely many left and right ideals, and every maximal subgroup of $T$ is a maximal subgroup of $S$ and thus it is defined by a finite complete rewriting system by assumption. Hence by induction $T$ is defined by a finite complete rewriting system. But now $S$ is an ideal extension of $T$ by the Rees quotient $S / T$ and $S/ T\cong J^0$ is a completely $0$-simple semigroup with finitely many left and right ideals, all of whose maximal subgroups are defined by finite complete rewriting systems. Thus by Theorem \ref{comp0simpExt}, $S/T$ is defined by a finite complete rewriting system. Therefore, applying Theorem~\ref{theoremIdealExt} we conclude that $S$ is defined by a finite complete rewriting system.
\end{proof}

In addition to this introduction, this paper comprises four sections. In
Section~2 we recall some basic definitions and results about string rewriting systems. 
Section~3 concerns complete rewriting systems for ideal extensions of semigroups, and is where we prove Theorem~\ref{theoremIdealExt}. 
Finally, the proof of Theorem~\ref{comp0simpExt} is given in Section~4 where complete rewriting systems for completely $0$-simple semigroups are considered.

\section{Preliminaries}
\label{sec_prelims}

In this section we give some  basic definitions and  results about string rewriting systems. 
For more background on string rewriting systems we refer the reader to \cite{BO93,TER03}.

Let $X $ be an alphabet. We denote by $X^*$ the free monoid on $X$
and by $X^+$ the free semigroup on $X$. For an element of $X^*$,
a {\it word} $w$, we denote the {\it length} of $w$ by $|w|$. Given a subset $Y$ of $X$, we write $|w|_Y$ for the total number of letters in $w$ that come from the subset $Y$.

A {\it presentation} is a pair $\langle X\mid R\rangle$, where  $X$ is an
alphabet and $R$ is a binary relation on $X^*$. The set $R$ is also referred to as rewriting system and its elements as {\it rewriting rules}. Usually, a rewriting
rule $r\in R $ is written in the form $r=(r_{+1},r_{-1})$ or,
simply, $r_{+1}\rightarrow r_{-1}$.  We say that a presentation $R$ is {\it finite} if both $R$ and $X$ are finite.

We define a binary relation $\ssr{R} $ on $X^*$, called a {\it
single-step reduction}, in the following way: $$u\ssr{R} v \
\Leftrightarrow \ u=w_1r_{+1} w_2 \ {\rm and} \ v=w_1r_{-1}w_2 $$
for some $(r_{+1},r_{-1}) \in R $ and $w_1,w_2\in X^* $. The
transitive and reflexive closure of $\ssr{R} $ is denoted by
$\RhCon{R}$. We use $\RhConAS{R}$ to denote the transitive closure of
$\ssr{R} $. 
A word $u\in X^*$ is said to be {\it $R$-reducible},
if there is a word $v\in X^*$ such that $u\ssr{R} v$. If a word is
not $R$-reducible, it is called {\it $R$-irreducible} or simply
irreducible. The set of all
 $R$-irreducible words is denoted by $Irr(R)$.

We denote by $\ThCon{R}$ the reflexive transitive symmetric closure
of
$\ssr{R}$ which is a congruence on the free monoid $X^*$, called the {\it Thue
congruence generated by $R$}.
The
quotient $X^* / \ThCon{R}$ of the free monoid $X^*$ by $\ThCon{R}$
is called the \emph{monoid defined
by $R$} and it is denoted by $M(X;R)$. The set $X$ is called the
\emph{generating set} and $R$ the set of \emph{defining relations}.  
A monoid is said to be defined by the presentation $\langle X\mid R\rangle$, or by the rewriting system
$R$, if $M\cong M(X;R)$. Thus, the elements of $M$ are identified
with congruence classes of words from $X^*$. Given two words $w,v \in X^*$ we write $w \equiv v$ if $w$ and $v$ are identical as words in $X^*$, and $w=v$ to mean they represent the same element of $M$ (that is, if $w / \ThCon{R} = v / \ThCon{R}$).

A rewriting system $R$ on $X $ is said to be {\it noetherian} if
the relation $\ssr{R} $ is well-founded, in other words, if there
are no infinite descending chains
$$w_1\ssr{R}w_2\ssr{R}w_3\ssr{R}\cdots
\ssr{R}w_n\ssr{R}\cdots . $$ It is called {\it confluent} if
whenever we have $u\RhCon{R}v$ and $u\RhCon{R}v'$ there is a word
$w\in X^*$ such that $v\RhCon{R}w$ and $v'\RhCon{R}w$. If $R$ is
simultaneously noetherian and confluent we say that $R$ is {\it
complete}. We say that a presentation is noetherian, confluent or complete if its associated rewriting system has the respective property.

It is easy to verify that, if $R$ is a noetherian rewriting
system, each congruence class  of $M(X;R)$  contains at least one
irreducible element. Assuming $R$ is noetherian, then $R$  is a
complete rewriting system if  and only if each congruence class of
$M(X;R)$ contains exactly one irreducible element; see \cite[Theorem 1.2.2]{TER03}.
Hence, a complete rewriting system fixes a unique normal form for
each of its congruence classes, given by taking the unique irreducible word in the class. 

There are obvious analogous definitions and results to those above obtained by replacing monoid by
semigroup, and the free monoid $X^*$ by the free semgiroup $X^+$, throughout. 
In this paper, we shall find it convenient to work in this slightly more general context of semigroups and semigroup presentations. 
Of course, all of the results obtained also hold for monoid presentations, since it is known, and quite straightforward to prove, that a monoid is defined by a finite complete monoid presentation if and only if it is defined by a finite complete semigroup presentation; see \cite{GM10} for details. All the main results here have been stated without reference to a specific type of presentation, and the aforementioned fact tells us that no ambiguity arises in doing so. 

It is important to note that for a finite rewriting system
$R$ on a set $X$, for each word $u$ of $X^+$ there are only finitely many single-step reductions that can by applied to $u$. 
If we also assume that $R$ is noetherian then $u$ has only finitely many descendants and so there is a maximum length that a descendant of $u$ can have. 
That maximum is called the
{\it stretch} of $u$ and it is denoted by $st_R(u)$. Observe that
if $u\rightarrow_R v$ then $st_R(u)\geq st_R(v)$. Also note that if $v$ is a proper factor of $u$ then $st_R(u)>st_R(v)$. 

We end this section with a technical result which will be used in the sequel. 
Let $S$ be a semigroup with a zero element (i.e. an element $b$ such that $bs = sb = b$ for all $s \in S$) defined by the complete semigroup
presentation $\langle X\mid R\rangle$. Suppose that the zero of $S$ is
represented by some word $z$ on $X^+$. Applying a Tietze
transformation we get a presentation $\langle
X\cup\{0\}\mid R\cup\{(z,0)\}\rangle$ which also defines the
semigroup $S$.
\begin{prop}\label{prop6.3}
Let $S$ be a semigroup with zero defined by the finite complete 
presentation $\langle X\mid R\rangle$. Let $0$ be a symbol not in
$X$ and let $z\in X^+$ be the irreducible element that represents
the zero. Then the rewriting system on $X\cup\{0\}$,
$$R_0 = R\cup\{(z,0)\}\cup\left\{(0x,0),(x0,0): x\in X\cup\{0\}\right\}$$ is finite
complete and defines $S$. Moreover, in this rewriting system $0$ is the irreducible representing the zero of $S$.
\end{prop}
\begin{proof}

First we will show that the rewriting system $R_0$ is  noetherian.
Suppose that there exists an infinite sequence $w_1 \ssr{{R_0}}
w_2 \ssr{{R_0}} \cdots \ssr{{R_0}} w_n \ssr{{R_0}}\cdots$. From
this sequence we 
construct another sequence $\widetilde{w}_1
\RhCon{{R}} \widetilde{w}_2 \RhCon{{R}} \cdots \RhCon{{R}}
\widetilde{w}_n \RhCon{{R}}\cdots$ in $\langle X\mid R\rangle$
in the following way. 

First, each of the words $\widetilde{w}_i$ in the new sequence is obtained from the corresponding $w_i$ in the original sequence by
replacing each occurrence of the symbol $0$ by $z$. Then each single step reduction from the original sequence is replaced by a reduction sequence as follows. 
For any $x\in X$,
since $z$ is $R$-irreducible, the single-step relation
$0x\ssr{{R_0}} 0$ has a corresponding non-empty sequence
$zx\RhConAS{R} z$ where the $0$ is replaced by the $z$. Thus
wherever in the original sequence a relation of the form $(0x,0)$ ($x\in X$)
is applied  we can replace it by a non-empty sequence with the $0$
replaced by $z$. An analogous statement holds for the relations of the form $(x0,0)$
and $(00,0)$. Also, if the relation $(z,0)$ is applied in the original sequence,
that is, if we have $w_i\ssr{{R_0}}w_{i+1}$, for some $i\in
\mathbb{N}_0$, where the relation used is $(z,0)$,  then 
$\widetilde{w}_i$ and $\widetilde{w}_{i+1}$ are identical as words and the single-step
reduction $w_i\ssr{{R_0}}w_{i+1}$ is replaced  by the empty
sequence. Otherwise, the single step relation used comes from $R$ (and so does not involve the letter $0$) in which case we simply apply the same rule in the new sequence. 

Now, $R$ is noetherian meaning that the sequence $$\widetilde{w}_1
\RhCon{{R}} \widetilde{w}_2 \RhCon{{R}} \cdots \RhCon{{R}}
\widetilde{w}_n \RhCon{{R}}\cdots$$ involves only a finite number of single-step reductions. 
Therefore, from the way that it was constructed we 
conclude that the number of  single-step reductions in the
sequence $$w_1 \ssr{{R_0}} w_2 \ssr{{R_0}} \cdots \ssr{{R_0}} w_n
\ssr{{R_0}}\cdots$$ that use relations from
$R\cup\left\{(0x,0),(x0,0): x\in X\cup\{0\}\right\}$  is finite
and hence those that use the relation $(z,0)$ is infinite. But this means that
for some $k\in \mathbb{N}_0$ all the single-step relations
$w_n\ssr{{R_0}}w_{n+1}$, for $n\geq k$, use the relation $(z,0)$. This is a contradiction because $z$ does not contain the letter $0$
and hence we can not apply the relation $(z,0)$ infinitely many
times to $w_k$. Therefore $R_0$ is noetherian.

Now, it is clear that the presentation $\langle X\cup\{0\}\mid
R_0\rangle$ defines the semigroup $S$ and that the
$R_0$-irreducible elements are $Irr(R)\backslash\{z\}\cup \{0\}$,
and hence they are in one-to-one correspondence with $S$. It
follows that $R_0$ is a finite complete
rewriting system defining $S$.
\end{proof}

\section{Complete rewriting systems for ideal extensions}

This section is dedicated to the proof of Theorem~\ref{theoremIdealExt}.

Let $S$ be a semigroup and let $T$ be an ideal of $S$. Let us
denote by $U$ the Rees quotient of $S$ by $T$. This means that $S$
is the ideal extension of $T$ by $U$. Suppose that $T$ and $U$
are defined by finite complete rewriting systems $\langle A\mid R\rangle$ and $\langle B\mid Q\rangle$, respectively.

By $B_0$
we denote the set  of generators in $B$
representing the zero of $U$. Note that $B_0$ could be empty. Also, denote by $Q_0$ the set of
rewriting rules $(u,v)$ in $Q$ such that $u$, and hence $v$, represents the zero
of $U$. Without loss of generality we may suppose that $B_0$ and $Q_0$ are non empty  and that $B_0$ contains a distinguished letter $0$ which is the unique irreducible word representing the zero of $U$. Indeed, if not we could use Proposition~\ref{prop6.3} to replace $\langle B\mid Q\rangle$ by a complete rewriting system also defining $U$ which has these properties.

For each word $u\in (B\backslash B_0)^+$, where $u$ represents the
zero in $U$, we fix a word $\rho(u)$ in $A^+$ such that the
relation $u=\rho(u)$ holds in $S$. Also, for each pair of letters
$a\in A$ and $b\in B\backslash B_0$, we fix words $\sigma(a,b)$
and $\pi(b,a)$ such that $ab=\sigma(a,b)$ and $ba=\pi(b,a)$
hold in $S$.

It is easy to see  that $S$ is defined
by the finite presentation ${\mathcal P}$  with generators $A\cup
B\backslash B_0$ and rewriting rules $R$, $Q\backslash Q_0$, and
\begin{eqnarray}
u& \rightarrow &\rho(u),  \label{rr1}\\
ab& \rightarrow &\sigma(a,b), \label{rr2}\\
ba & \rightarrow & \pi(b,a), \label{rr3}
\end{eqnarray}
where $a\in A$, $b\in B\backslash B_0$, $u\in (B\backslash
B_0)^+$ and either $(u,v)$ or $(v,u)$ is in $Q_0$.
For a proof see \cite[Proposition 4.4]{Rus99}.

We shall now show that moreover, from the assumption that  $\langle A\mid R\rangle$ and $\langle B\mid Q\rangle$ are
 complete, it follows that the presentation ${\mathcal P}$ is complete. Let us denote by $V$ the set of rewriting rules from $\mathcal P$.
 
\begin{lem}\label{lem1}
For any word $w\in (A\cup B\backslash B_0)^+$
\begin{enumerate}[(i)]
  \item if $w$ represents an element of $T$ then there exists $w'\in A^+$ such that $w\RhCon{V} w'$;
  \item otherwise $w\in (B\backslash B_0)^+$.
\end{enumerate}
\end{lem}
\begin{proof}
\noindent (i) Let $w$ be a word in $(A\cup B\backslash B_0)^+$ representing an element of $T$. 
If $w \in A^+$ then $w$ represents an element of $T$ and (i) holds trivially. 
Next, suppose that $w$ contains letters from both
alphabets $A$ and $B\backslash B_0$. It is clear that relations of the form (\ref{rr2}) and (\ref{rr3}) can be used to reduce $w$ to a word in the alphabet $A$.

Otherwise, $w$ is a word in $(B\backslash B_0)^+$.  Then since $w$ represents the zero of $U$ and because  $\langle B\mid Q\rangle$ has the form given in Proposition~\ref{prop6.3} there is a non empty reduction sequence $w\equiv w_0 \ssr{Q} w_1 \ssr{Q} \cdots \ssr{Q} w_n\equiv 0$. Let $k_0\in\{1,\ldots, n\}$ be the least $i$ such that $w_{i-1}\ssr{{Q_0}}w_i$. Since $w\in (B\backslash B_0)^+$ and by the choice of $k_0$ we conclude that the rewriting rule used in the reduction $w_{k_0-1}\ssr{{Q_0}}w_{k_0}$ must be of the form $u\rightarrow v$, where $(u,v) \in Q_0$. Hence the relation $u \rightarrow \rho(u)$ belongs to $\mathcal{P}$ and so we may replace $u\rightarrow v$ by $u \rightarrow \rho(u)$ to obtain a reduction sequence, in the rewriting system $V$, $w\equiv w_0 \ssr{{Q\backslash Q_0}}\ \ \cdots \ssr{{Q\backslash Q_0}}\ \ w_{k_0-1}
\ssr{{(1)}} \ z$, where $z$ has letters from  the alphabet $A$, and this puts us back in one of the cases considered in the previous paragraph.

\noindent (ii)  The statement follows from the fact that $T$ is an ideal of $S$, and hence if $w$ contains some letter from $A$ it would mean that $w$ represents an element of $T$.
\end{proof}

We use $\mathcal{M}(\mathbb{N})$ to denote the set of all finite multisets over the set of natural numbers $\mathbb{N}$. Recall that a multiset is like a set, but where are allows multiple occurrences of elements. We use the notation $[n_1,n_2,\ldots,n_r]$, where $n_i\in \mathbb{N}$ ($1\leq i\leq r$) to denote the elements of $\mathcal{M}(\mathbb{N})$, where permutations of the numbers in the list leave the element unchanged. 

The natural order $>$ on $\mathbb{N}$ induces an order $>_{mult}$ called the multiset order on $\mathcal{M}(\mathbb{N})$ 
where $M >_{mult} N$ if and only  if it is possible to transform $M$ into $N$ by carrying out the following procedure finitely many times: remove an element $x$ from $M$ and replace it by a finite number of elements all of which are smaller than $x$. For a more detailed explanation of multisets see \cite[Appendix A.6]{TER03}. It follows from \cite[Theorem A.6.5]{TER03} that $>_{mult}$ is a well-founded order on $\mathcal{M}(\mathbb{N})$.

\begin{lem}\label{Lemma7}
The relation  $\ssr{V}$ is noetherian.
\end{lem}
\begin{proof}
A word $w\in  (A\cup B\backslash B_0)^+$ can be written uniquely in the form $u_nv_n\cdots u_1v_1u_0$ with $n\in \mathbb{N}_0$, $u_0,u_n\in A^*$, $u_i\in A^+$, for $i=1,\ldots, n-1$, and $v_j\in  (B\backslash B_0)^+$, for $j=1, \ldots, n$.
Similarly a word $w'\in  (A\cup B\backslash B_0)^+$ can be written in the form $u'_mv'_m\cdots u'_1v'_1u'_0$.

We then write $w\prec w'$ if we have \begin{enumerate}[(i)]
  \item $[st_Q(v'_1), \ldots, st_Q(v'_m)] >_{mult}[st_Q(v_1), \ldots, st_Q(v_n)]$; or
  \item these multisets are equal and we have $v_0\equiv v'_0, \ldots, v_k\equiv v'_k$, for some $k<n=m$, and $v'_{k+1}\ssr{{Q \setminus Q_0}} \; \; \; v_{k+1}$; or
  \item the multisets are equal,  $v_0\equiv v'_0, \ldots, v_n\equiv v'_n$, and we have $u_0\equiv u'_0, \ldots, u_k\equiv u'_k$, for some $k<n=m$, and $u'_{k+1}\ssr{R}u_{k+1}$.
\end{enumerate}
Since  $>_{mult}$ is well-founded and both $\ssr{Q}$ and $\ssr{R}$ are noetherian it follows that the relation  $\prec$  is  a well-founded strict order on $(A\cup B\backslash B_0)^+$.

We claim that whenever $w'\ssr{V}w$ then $w\prec w'$ thus proving that the relation $\ssr{V}$ is noetherian. Consider the above decompositions for $w$ and $w'$. If $w'$ is reduced to $w$ by applying a relation from $R$ we will be in situation (iii). If a relation from $Q\backslash Q_0$ is applied then for some $j\in\{0, \ldots, n\}$ we have $v'_j\ssr{{Q\backslash Q_0}}\ \ v_j$ and hence $st_Q(v'_j)\geq st_Q(v_j)$, in which case, we are in situation (ii) if $st_Q(v'_j)= st_Q(v_j)$, and in situation (i) otherwise. 

Now suppose that when reducing $w'$ to $w$ a rule of the form $u\rightarrow \rho(u)$ is applied. Then,  for some $j\in\{0, \ldots, n\}$, a factor of $v'_j$ is replaced by $\rho(u)$. The letters of $v'_j$ which are not replaced appear in the word $w$. They give rise to  one or two proper factors of $v'_j$ if $u$ is a proper factor of $v'_j$, otherwise $v'_j$ simply disappears. In each of these cases the natural number $st_Q(v'_j)$ is replaced in the multiset  $[st_Q(v'_1), \ldots, st_Q(v'_m)]$ by a finite set of natural numbers all of which are smaller than $st_Q(v'_j)$. Thus situation (i) occurs.

The final possibility is that a rule of the form $ab\rightarrow \sigma(a,b)$ (or $ba\rightarrow \tau(b,a)$) is applied, which may be dealt with in a similar way to the previous case. Indeed, some factor $v'_j$ of $w'$ will be replaced by some suffix of it, and thus it has a smaller stretch. Hence  again situation (i) occurs. 
\end{proof}

\begin{proof}[Proof of Theorem~\ref{theoremIdealExt}]
As already observed, the presentation $\mathcal P$ defines the semigroup $S$. In Lemma~\ref{Lemma7} it was shown that the relation  $\ssr{V}$ is noetherian. Now, by Lemma~\ref{lem1} and the fact that presentations $\langle A\mid R\rangle$ and $\langle B\mid Q\rangle$ are
both complete, the irreducible elements with respect to $V$ are in one-to-one correspondence with the elements of $S$. Consequently, $\mathcal P$ is a finite complete presentation defining $S$.
\end{proof} 

\section{Complete rewriting systems for completely 0-simple semigroups}
\label{sec_comp0simple}
This section will be devoted to the proof of Theorem~\ref{comp0simpExt}.  
We recall that by the Rees-Suschkevitz Theorem \cite[Section~3.2]{How95} a completely $0$-simple semigroup $S$ is isomorphic to a \emph{$0$-Rees matrix semigroup} $M^0[G,I,\Lambda,P]$, where $G$ is a group isomorphic to any (and hence all) non-zero maximal subgroups of $S$, $I$ is a set indexed by the set of all $0$-minimal right ideals of $S$, $\Lambda$ is a set indexed by the set of all $0$-minimal left ideals of $S$, and $P = (p_{\lambda i})$ is a regular $\Lambda \times I$ matrix with entries from $G \cup \{ 0 \}$. Multiplication in $M^0[G,I,\Lambda,P] = (I \times G \times \Lambda) \cup \{ 0 \}$ is given by
\[
\begin{array}{c}
(i,g,\lambda)(j,h,\mu) =
\begin{cases}
(i, gp_{\lambda j}h, \mu) & \mbox{if $p_{\lambda j} \neq 0$} \\
0 &  \mbox{if $p_{\lambda j} = 0$} \\
\end{cases} \\
0 (i,g,\lambda) = (i,g,\lambda)0 = 00 = 0.
\end{array}
\]
There is an analogous construction for completely simple semigroups, given by taking Rees matrix semigroups $M[G,I,\Lambda,P] = I \times G \times \Lambda$ over groups, where $P$ is a matrix with entries from $G$, and multiplication is given by $(i,g,\lambda)(j,h,\mu) = (i, gp_{\lambda j}h, \mu)$.

Let $S = M^0[G,I,\Lambda,P]$ be an arbitrary completely $0$-simple semigroup. Without loss of generality we may suppose that $1 \in I$, $1 \in \Lambda$ and that $p_{1 1} \neq 0$ so that $\{ 1 \} \times G \times \{ 1 \}$ is a group $\gh$-class. Let $\lb A | R \rb$ be a complete semigroup presentation for the group $G$. Let $e \in A^+$ be a fixed word representing the identity element of $G$. Let $B = \{ b_i : i \in I \setminus \{ 1 \} \}$, $C = \{ c_{\lambda} : \lambda \in \Lambda \setminus \{ 1 \}  \}$ and $X = A \cup B \cup C \cup \{ 0 \}$. 

In \cite[Theorem~6.2]{HR94} it is proved that $\langle A, B, C \;|\; R, \eqref{rr4}, \eqref{rr5} \rangle$ (where \eqref{rr4} and \eqref{rr5} are defined below) is a presentation for $S$ as a semigroup with zero. This may be converted into a genuine semigroup presentation for $S$ by adding a new generating symbol $0$ and relations \eqref{rr6}. Adding to this two additional families of redundant relations \eqref{rr5'} and \eqref{rr5''}, we obtain the following presentation for $S$: 
\begin{eqnarray}
\nonumber  \mathcal{P}_S = \langle X\ | \ R,  & &    \\
b_i e \rightarrow b_i, \ ec_{\lambda} \rightarrow c_{\lambda} & & (i \in I\backslash\{1\}, \lambda \in \Lambda\backslash\{1\}) \label{rr4}  \\
eb_i \rightarrow p_{1 i}, \
c_{\lambda} e \rightarrow p_{\lambda 1}, \ c_{\lambda} b_i \rightarrow p_{\lambda i}
 & & (i \in I\backslash\{1\}, \lambda \in \Lambda\backslash\{1\}) \label{rr5}\\
ab_i \rightarrow ap_{1 i}, \ c_\lambda a \rightarrow p_{\lambda 1} a  
& &
(a \in A, i  \in I \backslash \{ 1 \}, \lambda \in \Lambda \backslash \{ 1 \})
\label{rr5'} \\
c_{\lambda} c_{\mu} \rightarrow p_{\lambda 1} c_\mu, \ b_i b_j \rightarrow b_i p_{i j} & & 
(i, j  \in I \backslash \{ 1 \}, \lambda, \mu \in \Lambda \backslash \{ 1 \})
\label{rr5''} \\ 
x0 \rightarrow 0, \ 0x\rightarrow 0 & & (x\in X) \label{rr6}\rangle
\end{eqnarray}
Note the slight abuse of notation in this presentation: the symbols $p_{\lambda i}$ appearing in the relations are really fixed words from $A^+\cup\{0\}$ representing the elements $p_{\lambda i }$ of $G\cup\{0\}$. We will denote the rewriting rules of this presentation by $U$.

Our aim is to prove that $\mathcal{P}_S$ is in fact complete. We first need some technical lemmas about the presentation. 

\begin{lem}\label{dagger}
\begin{enumerate}[(i)]
\item The word $b_i w c_\lambda$ with $i \in I \setminus \{ 1 \}$, $\lambda \in \Lambda \setminus \{ 1 \}$ and $w \in A^*$ represents the element $(i,g,\lambda)$ of $S$, where $g \in G$ is the element represented by $w$ if $w \in A^+$, or $g=1$ if $w$ is the empty word.
\item The word $b_i w$ with $i \in I \setminus \{ 1 \}$ and $w \in A^*$ represents the element $(i,g,1)$ of $S$, where $g \in G$ is the element represented by $w$ if $w \in A^+$, or $g=1$ if $w$ is the empty word.
\item The word $w c_\lambda$ with $\lambda \in \Lambda \setminus \{ 1 \}$ and $w \in A^*$ represents the element $(1,g,\lambda)$ of $S$, where $g \in G$ is the element represented by $w$ if $w \in A^+$, or $g=1$ if $w$ is the empty word.
\item The word $w \in A^+$ represents the element $(1,g,1)$ where $g \in G$ is the element represented by the word $w$.
\end{enumerate}
\end{lem}
\begin{proof}
This follows from the construction of the presentation $\mathcal{P}_S$, together with the proof of \cite[Theorem~6.2]{HR94}.
\end{proof}

The next lemma identifies the normal forms of ${\mathcal P}_S$. First some notation. We write $B^1 A^* C^1$ to denote the set of all non-empty words from the set
\[
\{ bwc : b \in B \cup \{ 1 \}, \ c \in C \cup \{ 1 \} \ \& \ w \in A^*  \}.
\]
Likewise we use the notation $B^1 A^*$ and $A^* C^1$.
\begin{lem}\label{Lemma1}
Let $w \in X^+$ be arbitrary.
\begin{enumerate}[(i)]
\item If $w$ represents an element of $S \backslash \{0\}$ then there is a word  $w' \in B^1 A^* C^1$ such that $w\RhCon{U} w'$. 
\item If $w$ represents the zero of $S$ then $w\RhCon{U} 0$.
\end{enumerate}
\end{lem}
\begin{proof}
We prove (i) and (ii) simultaneously. The proof is by induction on the total number $|w|_{B \cup C}$ of letters in the word that come from the set $B \cup C$. If $|w|_{B \cup C}=0$ then either $w$ contains the letter $0$, in which case there is an obvious reduction from $w$ to $0$, or $w \in A^+$, so  $w$ represents an element of $G$ and is already written in the required form, and we are done by setting $w' \equiv w$.

Now suppose that $|w|_{B \cup C} > 0$. If $w \in B^1 A^* C^1 \cup \{0\}$ then we are done by setting $w' \equiv w$  and using the fact that the words in $B^1 A^* C^1$ all represent nonzero elements of $S$. Next suppose that $w \not\in B^1 A^* C^1 \cup \{0\}$. If $w$ contains the letter $0$ there is an obvious reduction to $0$. Otherwise, $w$ must have a subword which is the left hand side of one of the rewriting rules \eqref{rr5'}, \eqref{rr5''} or the third relation in \eqref{rr5}, and applying this rule will result in a word $v$ satisfying $|v|_{B \cup C} < |w|_{B \cup C}$, and the result follows by induction. 
\end{proof}

\begin{lem}
The relation $\ssr{U}$ is noetherian.
\end{lem}
\begin{proof}
An arbitrary  word $w$ in $X^+$ has the form $x_nu_nx_{n-1}\cdots
x_1u_1x_0$, with $n\in\mathbb{N}_0$, $u_i\in B\cup C\cup \{0\}$, and $x_j\in
A^*$.

Following the approach as in Lemma~\ref{Lemma7} we will define a well-founded strict relation $\prec$ on the set $X^+$ such that $w'\ssr{U} w$ implies $w\prec w'$ and hence prove that $\ssr{U}$ is noetherian. 
Let $w,w'\in X^+$ and decompose them as $w\equiv x_nu_nx_{n-1}\cdots
x_1u_1x_0$ and $w'\equiv x'_mu'_mx'_{m-1}\cdots x'_1u'_1x'_0$.
We say that $w \prec w'$ if:
\begin{enumerate}[(i)]
\item $|w'|_{B\cup C}> |w|_{B\cup C}$; or
\item if $|w'|_{B\cup C}= |w|_{B\cup C}$ and $|w'|_{\{0\}}> |w|_{\{0\}}$; or
\item  if $|w'|_{B\cup C}= |w|_{B\cup C}$, $|w'|_{\{0\}}= |w|_{\{0\}}$ (thus $n=m$) and 
\[
[st_R(x'_1), \ldots, st_R(x'_m)] >_{mult}[st_R(x_1), \ldots, st_R(x_n)],
\]
where $>_{mult}$ is the multiset order on $\mathcal{M}(\mathbb{N})$; or
  \item if  $|w'|_{B\cup C}= |w|_{B\cup C}$, $|w'|_{\{0\}}= |w|_{\{0\}}$, the multisets are equal and we have $x_0\equiv x'_0, \ldots x_k\equiv x'_k$, for some $k<n=m$, and $x'_{k+1}\ssr{R}x_{k+1}$.
\end{enumerate}
The natural order $>$ on $\mathbb{N}$ and the multiset order $>_{mult}$ on $\mathcal{M}(\mathbb{N})$ are both well-founded and the relation $\ssr{R}$ is noetherian since $\lb A | R \rb$ is a complete presentation. It then readily follows that the above relation $\prec$ on the set $X^*$ is a well-founded strict order. We will now show that whenever $w'\ssr{U} w$ we get $w\prec w'$ and hence prove that $\ssr{U}$ is noetherian.

Suppose that the $w'$ is reduced to $w$ in one step. If a rewriting rule from $R$ is applied then either situation (iii) or (iv) occurs thus showing that $w\prec w'$. If a rewriting rule of type (\ref{rr4}) is applied then situation (iii) occurs since for some $i\in \{0,\ldots,m\}$, $x_i$ is a proper factor of $x_i'$. The easiest case is when a rewriting rule of type (\ref{rr5}), \eqref{rr5'} or \eqref{rr5''} is applied because situation (i) occurs.  We distinguish three cases when a rewriting rule of type (\ref{rr6}) is applied: if $x\in B\cup C$ we are in situation (i); if $x\equiv 0$ then situation (ii) occurs; otherwise, if $x\in A$ we have (iii).

This covers all possible rewrite rules and so completes the proof of the lemma. 
\end{proof}

\begin{proof}[Proof of Theorem~\ref{comp0simpExt}]
Suppose that $G$ is defined by a finite complete rewriting system $\langle A\mid R\rangle$. Then the semigroup $S=M^0[G;I,\Lambda;P]$  is defined by the presentation ${\mathcal P}_S=\langle X\mid U\rangle$. In the previous lemma we have seen that the relation $\ssr{U}$ over $X$ is noetherian.

We claim that the irreducible elements of $U$ are in one-to-one correspondence with the elements of $S$, and thus  ${\mathcal P}_S$ is a finite complete presentation defining $S$. Indeed, let $\alpha, \beta\in X^+$ be irreducible with respect to $\ssr{U}$ and suppose that $\alpha = \beta$ in $S$. We must show that $\alpha\equiv \beta$. If $\alpha$ (and hence $\beta$) represents zero then, since $\alpha$ and $\beta$ are both irreducible, Lemma~\ref{Lemma1}(ii) implies $\alpha\equiv \beta\equiv 0$. Otherwise, since $\alpha$ and $\beta$ are irreducible, Lemma~\ref{Lemma1} implies that $\alpha\equiv b\gamma c$ and $\beta\equiv b' \gamma' c'$, for some $b,b'\in B\cup\{1\}$, $c,c'\in C\cup\{1\}$ and $\gamma, \gamma'\in A^*$. Let $(i,g, \lambda)$ be the element of $S$ represented by $\alpha$ (and hence also $\beta$ since $\alpha = \beta$). Since $b\gamma c\equiv \alpha =\beta\equiv b' \gamma' c'$ it follows by Lemma~\ref{dagger} that $b\equiv b'$ and $c\equiv c'$.

If $b\equiv c\equiv 1$ then by Lemma~\ref{dagger}, $\gamma$ and $\gamma'$  are non-empty words in $G$ representing the same element.
 Since $b\gamma c$ is irreducible with respect to $\ssr{U}$ and $R\subseteq U$, it follows that $\gamma$ must be irreducible with respect to $\ssr{R}$. Likewise, $\gamma'$ is irreducible with respect to $\ssr{R}$. But by assumption  $\langle A\mid R\rangle$ is complete, therefore $\gamma\equiv \gamma'$ and thus $\alpha\equiv \beta$ as required.

If  $b$ or $c$ are non-empty we can distinguish two cases depending of whether $g$ is the identity of $G$ or not. If $g$ is the identity we conclude that $\gamma$ and $\gamma'$ are empty words. Indeed, if $\gamma$ is not the empty word we know by Lemma~\ref{dagger} that $\gamma$ represents the identity and hence, since   $b\gamma c$ is irreducible with respect to $\ssr{U}$, and in particular with respect to $\ssr{R}$,
we get $\gamma\equiv e$. But then we get a contradiction since by relations (\ref{rr4}) we have $bec\ssr{U} bc$. Thus $\gamma$ and analogously $\gamma'$  are empty which means that $\alpha\equiv \beta$.

Now if $b$ or $c$ are non-empty and $g$ is not the identity we conclude by Lemma~\ref{dagger} that $\gamma$ and $\gamma'$ are both non-empty. As before, we conclude that $\gamma$ and $\gamma'$ are irreducible with respect to $\ssr{R}$ and therefore $\gamma\equiv \gamma'$  which proves that $\alpha\equiv \beta$ as required.
\end{proof} 

\bibliographystyle{abbrv}

\end{document}